\documentclass[twoside,leqno]{amsart}

\usepackage{amsmath,color,amsthm}

\numberwithin{equation}{section}

\newtheorem{theorem}{Theorem}[section]
\newtheorem{corollary}[theorem]{Corollary}
\newtheorem{lemma}[theorem]{Lemma}
\newtheorem{proposition}[theorem]{Proposition}

\newtheorem{example}[theorem]{\sl Example}
\newtheorem{definition}[theorem]{\sl Definition}

\theoremstyle{definition}
\newtheorem{Remark}[theorem]{Remark}

\newcommand{\beqn}{\begin{eqnarray}}
\newcommand{\eeqn}{\end{eqnarray}}
\newcommand{\beqnn}{\begin{eqnarray*}}
\newcommand{\eeqnn}{\end{eqnarray*}}

\newcommand{\lf}{\left\lfloor}

\newcommand{\rf}{\right\rfloor}

\newcommand{\EE}{{\bf  E}}

\newcommand{\PP}{{\bf  P}}

\newcommand{\Var}{{\bf Var}}
\newcommand{\Cov}{{\bf Cov}}

\newcommand{\Lc}{{\mathcal L}}
\newcommand{\Lto}{{\stackrel{\Lc}{\to}}}
\newcommand{\Leq}{{\,\stackrel{\Lc}{=}\,}}

\newcommand{\begp}{\begin{proposition}}
\newcommand{\enp}{\end{proposition}}
\newcommand{\begt}{\begin{theorem}}
\newcommand{\ent}{\end{theorem}}
\newcommand{\begl}{\begin{lemma}}
\newcommand{\enl}{\end{lemma}}
\newcommand{\begc}{\begin{corollary}}
\newcommand{\enc}{\end{corollary}}
\newcommand{\begcl}{\begin{claim}}
\newcommand{\encl}{\end{claim}}
\newcommand{\begr}{\begin{Remark}\rm}
\newcommand{\enr}{\end{Remark}}
\newcommand{\begal}{\begin{algorithm}}
\newcommand{\enal}{\end{algorithm}}
\newcommand{\begd}{\begin{definition}}
\newcommand{\enf}{\end{definition}}
\newcommand{\begx}{\begin{example}}
\newcommand{\enx}{\end{example}}
\newcommand{\bega}{\begin{array}}
\newcommand{\ena}{\end{array}}

\newcommand{\sfrac}[2]{{\textstyle\frac{#1}{#2}}}

\def\rompar(#1){\textup(#1\textup)}    

\newcommand{\refS}[1]{Section~\ref{#1}}
\newcommand{\refT}[1]{Theorem~\ref{#1}}

\newcommand{\refL}[1]{Lemma~\ref{#1}}
\newcommand{\refP}[1]{Proposition~\ref{#1}}

\newcommand\ie{i.e.\spacefactor=1000}

\newcommand\QuickSort{\texttt{QuickSort}}

\newcommand{\ignore}[1]{}

\begin{document}

\newcommand{\tab}[0]{\hspace{.1in}}

\title[Limit Distn.\ for QuickSort Symbol Comparisons is Nondegenerate]
{The Limiting Distribution for the Number of Symbol Comparisons Used by QuickSort is Nondegenerate (Extended Abstract)}

\author{\ \ \ \ Patrick Bindjeme\ \ \ \ \vspace{.2cm}\\ James Allen Fill}
\address{Department of Applied Mathematics and Statistics,
The Johns Hopkins University,
34th and 
Charles Streets,
Baltimore, MD 21218-2682 USA}
\email{bindjeme@ams.jhu.edu and jimfill@jhu.edu}
\thanks{Research supported by the Acheson~J.~Duncan Fund for the Advancement of Research in Statistics.}

\date{January~27, 2012}

\maketitle

\vspace{-.2in}
\begin{center}
{\sc Abstract}
\vspace{.2cm}
\end{center}

\begin{small}
In a continuous-time setting, Fill~\cite{fjfeb2010} proved, for a large class of probabilistic sources, that the number of symbol comparisons used by {\tt QuickSort}, when centered by subtracting the mean and scaled by dividing by time, has a limiting distribution, but proved little about that limiting random variable~$Y$---not even that it is nondegenerate.  We establish the nondegeneracy of~$Y$.  The proof is perhaps surprisingly difficult.
\end{small}
\vspace{-.03in}

\section{The number of symbol comparisons used by {\tt QuickSort}:\ Brief review of a limiting-distribution result}
\label{revng}
In this section we briefly review the main theorem of~\cite{fjfeb2010}.  An infinite sequence of independent and identically distributed keys is generated; each key is a random word 
$(w_1, w_2, \ldots) = w_1 w_2 \cdots$, that is, an infinite sequence, or ``string'', of symbols $w_i$ drawn from a totally ordered finite alphabet~$\Sigma$.  The common distribution~$\mu$ of the keys (called a \emph{probabilistic source}) is allowed to be any distribution over words, \ie,\ the distribution of any stochastic process with time parameter set $\{1, 2, \dots\}$ and state space~$\Sigma$ .  
We know thanks to Kolmogorov's consistency criterion (e.g.,\ Theorem~$3.3.6$ in \cite{c2001}) that the possible distributions~$\mu$ are in one-to-one correspondence with consistent specifications of finite-dimensional marginals, i.e.,\ of the \emph{fundamental probabilities}
\begin{equation}
\label{defini}     
p_w := \mu(\{w_1w_2\cdots w_k\}\times \Sigma^\infty)\,\,\text{with}\,\, w = w_1w_2\cdots w_k \in \Sigma^*.
\end{equation}
This $p_w$ is the probability that a word drawn from~$\mu$ has~$w$ as its length-$k$ prefix.

For each~$n$, Hoare's~\cite{h1962} \QuickSort\ algorithm can be used to sort the first~$n$ keys to be generated.  We may and do assume that the \emph{first} key in the sequence is chosen as the pivot, and that the same is true recursively (in the sense, for example, that the pivot used to sort the keys smaller than the original pivot is the \emph{first} key to be generated that is smaller than the original pivot).  A comparison of two keys is done by scanning the two words from left to right, comparing the symbols of matching index one by one until a difference is found.
We let $S_n$ denote the total number of symbol comparisons needed when~$n$
keys are sorted by {\tt QuickSort}. 

\begt[Fill~\cite{fjfeb2010}, Theorem~3.1]
\label{T:main}
Consider the continuous-time setting in which keys are generated from a probabilistic source at the arrival times of an independent Poisson process~$N$ with unit rate.  Let $S(t) = S_{N(t)}$ denote the number of symbol comparisons required by {\tt QuickSort} to sort the keys generated through epoch~$t$, and let
\begin{equation}
\label{yt}
Y(t) := \frac{S(t) - \EE\,S(t)}{t}, \qquad 0 < t < \infty.
\end{equation}
Assume that
\begin{equation}
\label{cond}
\sum_{k = 0}^{\infty} \Big( \sum_{w \in \Sigma^k} p^2_w \Big)^{1 / 2} < \infty
\end{equation}
with $p_w$ as in~\eqref{defini}.
Then there exists a random variable $Y$ such that $Y(t) \to Y$ in $L^2$ as $t \to \infty$.  In particular, $Y(t)\, \Lto \,Y$ and {\rm [}because $\EE\,Y(t) \to \EE\,Y${\rm ]} we have $\EE\,Y = 0$, and $\Var\,Y(t) \to \Var\,Y$.
\ent

In the full-length paper (in preparation) corresponding to~\cite{fjfeb2010}, this theorem will be extended by replacing the power $1 / 2$ in~\eqref{cond} by $1 / p$ for any given $p \in [2, \infty)$ and concluding that $Y(t) \to Y$ in $L^p$.

From \refT{T:main} we know that $\Var\,S(t) = O(t^2)$ as $t \to \infty$, but we don't know that $\Var\,S(t) = \Theta(t^2)$ because the theorem does not contain the important information that the limiting random variable~$Y$ is nondegenerate (\ie,\ does not almost surely vanish).  {\bf The purpose of the present extended abstract is to show 
that~$Y$ is nondegenerate}; this is stated as our main \refT{mainngthm} below.  The proof turns out to be surprisingly difficult; we do not know the value of $\Var\,Y$, and the proof of \refT{mainngthm} does not provide it.  The consequence $\Var\,S(t) = \Theta(t^2)$ of our \refT{mainngthm} settles a question that has been open since the work of Fill and 
Janson~\cite{fjbits2004} even in the special case of the standard binary source with $\Sigma = \{0, 1\}$ and the fundamental probabilities of~\eqref{defini} equal to $2^{-k}$.

\section{Main results}\label{mainng}

The following is the main theorem of this extended abstract.

\begt\label{mainngthm}
The limit distribution in Theorem \ref{T:main} is nondegenerate.
\ent

Throughout this extended abstract, we work in the setting of \refT{T:main}.  \refT{mainngthm} follows immediately from Propositions \ref{mainngthm1}--\ref{mainngthm2} in this section.

\begin{definition}
\label{Skwngdf}
\emph{
For an integer~$k$ and a prefix $w \in \Sigma^k$ we define (with little possibility of notational confusion),
for comparisons among keys that have arrived by epoch~$t$, the counts
\begin{align*}
S_k(t)
 &:= \mbox{number of comparisons of $(k + 1)$st symbols}, \\
S_w(t)
 &:= \text{number of comparisons of $(k + 1)$st symbols
between keys with prefix~$w$}.
\end{align*}
}
\end{definition}

The following two propositions combine to establish \refT{mainngthm}.  We write 
$\Sigma^* := \cup_{0 \leq k < \infty} \Sigma^k$ for the set of all prefixes.

\begp\label{mainngthm1}
If the random variables $S_w(t)$, $w \in \Sigma^*$, are nonnegatively correlated for each fixed $t$, then the limit distribution in Theorem \ref{T:main} is nondegenerate.
\enp

A proof of \refP{mainngthm1} can be found in \refS{needed} (see Subsection~\ref{condng}).

\begp\label{mainngthm2}
For each fixed $t$, the random variables $S_w(t)$, $w \in \Sigma^*$, are nonnegatively correlated.
\enp

Following six lemmas in \refS{S:lemmas}, a proof of \refP{mainngthm2} can be found in \refS{mainr} (specifically:\ in Subsection~\ref{g1}).
%
%
%

\section{Proof of Proposition \ref{mainngthm1}}
\label{needed}

\subsection{A lower bound for the variance of $K(t)$}
\label{Ktvarng}

\begin{definition}
\label{defkt}
\emph{
If $K_n$ is the number of key comparisons needed to sort the first~$n$ keys to arrive using {\tt Quicksort}, and~$N$ is the Poisson process in \refT{T:main} (independent of the generation of the keys), we define $K(t) := K_{N(t)}$.
}
\end{definition}

In order to prove Proposition \ref{mainngthm1}, we first establish the following lemma.

\begl\label{Kbdng} We have
\begin{eqnarray*}
 \Var\,K(t) \ge (1+ o(1))\sigma^2t^2\,\,\text{as}\,\, t \to \infty
\end{eqnarray*}
where $\sigma^2 := 7 - \frac{2}{3}\pi^2$.
\enl

\begin{proof}
By the law of total variance (namely, variance equals the sum of expectation of conditional variance and variance of conditional expectation) we have
\begin{equation}
\label{lawtotvar}
\Var\,K(t) \geq \EE\,\Var[K(t)\,\vert\, N(t)] = e^{-t}\sum_{n=0}^{\infty}\frac{t^n}{n!} \Var\,K_n.
\end{equation}
From (for example) $(1.2)$ in~\cite{fj2002} we have
$$
\Var\,K_n = 7 n^2 - 4 (n + 1)^2 H_n^{(2)} - 2 (n + 1) H_n +13 n,
$$
so
$$
\lim_{n \to \infty} \frac{\Var\,K_n}{n^2} = \sigma^2.
$$
It follows that, given $\alpha > 0$, there exists $n_\alpha$ such that
$$
\Var\,K_n \geq (1 - \alpha) \sigma^2 n^2\,\,\text{for all}\,\,n \geq n_\alpha.
$$
We therefore have from~\eqref{lawtotvar} that
\begin{eqnarray*}
\Var\,K(t) & \ge & (1-\alpha)\sigma^2 e^{-t}\sum_{n=n_\alpha}^{\infty}\frac{t^n}{n!}n^2 \\
& = & (1+o(1))(1-\alpha)\sigma^2 e^{-t}\sum_{n=0}^{\infty}\frac{t^n}{n!}n^2\,\,\text{as}\,\, t \to \infty \\
& = & (1+o(1))(1-\alpha)\sigma^2(t^2+t)
= (1+o(1))(1-\alpha)\sigma^2t^2.
\end{eqnarray*}
Since $\alpha > 0$ is arbitrary, the lemma follows.
\end{proof}

\subsection{Proof of Proposition \ref{mainngthm1}}\label{condng}

\begin{definition}
\label{defykyty}
\emph{
For any nonnegative integer $k$, with $S_k(t)$ as in Definition \ref{Skwngdf} we define
\begin{eqnarray*}
Y_k(t) & := & \frac{S_k(t) - \EE\,S_k(t)}{t}.
\end{eqnarray*}
}
\end{definition}

\begin{proof}[\textbf{Proof of Proposition \ref{mainngthm1}}] With $Y(t)$ as in \refT{T:main}, we have
\begin{eqnarray*}
Y(t) & = & \sum_{k=0}^\infty \, Y_k(t)
\end{eqnarray*}
and, from \refT{T:main},
\begin{eqnarray}
\Var\,Y(t)\,\to\,\Var\,Y\,\, \text{as}\,\, t\,\to\,\infty. \label{equa4}
 \end{eqnarray}

Knowing that
$$
\EE\,Y_k(t) = 0\,\,\text{for any nonnegative integer}\,\, k\,\,\text{and}\,\,t \in (0,\infty),
$$
that
$$
\EE\,Y(t) = 0\,\,\text{for}\,\,t \in (0,\infty),
$$
and finally that the random variables $Y_k(t)$ satisfy the hypotheses of the elementary probabilistic Lemma~2.8 
of~\cite{fjfeb2010} for $p_0 = 2$, we have for any $t \in (0,\infty)$ that
\begin{equation}
\label{lim1}
\Var\,\left( \sum_{k=0}^n Y_k(t) \right) \to \Var\,Y(t)\,\,\text{as}\,\, n \to \infty.
\end{equation}
Now, from the fact that
$$
S_k(t) = \sum_{w \in \Sigma^k} \, S_w(t)\,\,\text{for any}\,\,k,
$$
we have
\begin{eqnarray*}
\Var \left( \sum_{k=0}^n Y_k(t) \right) 
& = & \sum_{k=0}^n \Var\,Y_k(t) + 2 \sum_{0 \leq i < j \leq n} \Cov(Y_i(t), Y_j(t)) \\
& = & \sum_{k=0}^n \Var\,Y_k(t) + \frac{2}{t^2} \sum_{0 \leq i < j \leq n} \Cov(S_i(t), S_j(t)) \\
& = & \sum_{k=0}^n \Var\,Y_k(t) + \frac{2}{t^2} \sum_{0 \leq i < j \leq n} \sum_{\genfrac{}{}{0pt}{}{w \in \Sigma^i}{w' \in \Sigma^j}} \Cov(S_w(t), S_{w'}(t)).
\end{eqnarray*}
This allows us to conclude that if for each fixed~$t$ the random variables $S_w(t)$ with $w \in \Sigma^*$ are nonnegatively correlated, then
$$
\Var \left( \sum_{k=0}^n Y_k(t) \right) \geq \sum_{k=0}^n \Var\,Y_k(t)
$$
and therefore, considering~\eqref{lim1}, that
\begin{equation}
\label{varlim}
\Var\,Y(t) \geq \sum_{k=0}^\infty \Var\,Y_k(t) \geq \Var\,Y_0(t).
\end{equation}
As noted in~\cite[(3.3)--(3.4)]{fjfeb2010}, 
for any fixed $t$ and any $k \geq 0$ we have that
\begin{equation}
\label{ind}
\mbox{the r.~variables $S_w(t)$ with $w \in \Sigma^k$ are independent, and $S_w(t) \Leq K(p_wt)$},
\end{equation}
where $K(t)$ is defined in Definition~\ref{defkt}. It follows from~\eqref{varlim} and \refL{Kbdng} that
$$
\Var\,Y(t) \geq t^{-2} \Var\,K(t) \geq (1 + o(1)) \sigma^2,
$$
which implies from~\eqref{equa4} that
$$
\Var\,Y \geq \sigma^2 > 0.
$$
\end{proof}

\section{Six lemmas} 
\label{S:lemmas}

In the next section, we will need the following six lemmas.  We write $\kappa_n := \EE\,K_n$ for the expected number of key comparisons required to sort the first~$n$ keys to arrive.

\begl
\label{delt2lmng}
If
$$
\Delta_2(n, a, b) := 2 \sum_{j=1}^{b} \frac{1}{n + 2 - j - a}(b - 1 + \kappa_{j - 1} + \kappa_{b - j} - \kappa_b)
$$
for any nonnegative integers $a$, $b$, and $n$ with $a+b \le n$, then we have
\begin{eqnarray*}
\Delta_2(n,a,b) \ge 0.
\end{eqnarray*}
\enl

\begin{proof}
It is well known, and can easily be checked, using the explicit formula $\kappa_n = 2 (n + 1) H_n - 4 n$, that
\begin{align}
\label{orig1}
& b - 1 + \kappa_{j - 1} + \kappa_{b - j} - \kappa_b\mbox{\ is symmetric about $j = (b + 1) / 2$}, \\
 & \quad\qquad\text{and decreasing in $j = 1, \ldots, \lf (b + 1) / 2 \rf$}\nonumber.
\end{align} 
Thus
\begin{eqnarray}
\lefteqn{\Delta_2(n,a,b)} \nonumber \\
 & = & 2\sum_{j=1}^{\lfloor b/2 \rfloor}\left[\frac{1}{n+2-a-j} + \frac{1}{n+2-a-(b+1-j)}\right] \label{delta2}\\
 & & \mbox{}\qquad\qquad\times  (b - 1 + \kappa_{j-1} + \kappa_{b - j} - \kappa_b)\nonumber\\
 &  & \mbox{}\qquad+\mathbf{1}(b\,\text{is odd}) \frac{2}{n+2-a-(b+1)/2}(b - 1 + 2\kappa_{(b-1)/2} - \kappa_b)\nonumber.
 \end{eqnarray}
The first of the two terms in~\eqref{delta2} equals
\begin{eqnarray}
\label{delta21}
\lefteqn{\qquad 2[2(n+2-a)-(b+1)]} \\
& & \times \sum_{j=1}^{\lfloor b/2 \rfloor} \frac{1}{(n+2-a-j)(n+2-a-(b+1-j))} 
(b - 1 + \kappa_{j-1} + \kappa_{b - j} - \kappa_b).\nonumber
\end{eqnarray}
Each of the two factors in the sum in~\eqref{delta21}, namely $\frac{1}{(n+2-a-j)(n+2-a-(b+1-j))}$ and $b - 1 + \kappa_{j-1} + \kappa_{b - j} - \kappa_b$, decreases in $j$ over the range of summation, so by ``Chebyshev's other inequality" \cite{MR789242} we have
\begin{eqnarray*}
\lefteqn{\hspace{-.5in}\frac{1}{\lfloor b/2 \rfloor}\sum_{j=1}^{\lfloor b/2 \rfloor}\frac{1}{(n+2-a-j)(n+2-a-(b+1-j))}
(b - 1 + \kappa_{j-1} + \kappa_{b - j} - \kappa_b)} \\
& \ge & \left[\frac{1}{\lfloor b/2 \rfloor}\sum_{j=1}^{\lfloor b/2 \rfloor}\frac{1}{(n+2-a-j)(n+2-a-(b+1-j))}\right]\\
& & \mbox{}\times  \left[\frac{1}{\lfloor b/2 \rfloor}\sum_{j=1}^{\lfloor b/2 \rfloor}(b - 1 + \kappa_{j-1} + \kappa_{b - j} - \kappa_b)\right].
\end{eqnarray*}
Thus~\eqref{delta21} is greater than or equal to
\begin{eqnarray*}
\lefteqn{2[2(n+2-a)-(b+1)]\lfloor b/2 \rfloor}\\
& & \mbox{}\times\left[\frac{1}{\lfloor b/2 \rfloor}\sum_{j=1}^{\lfloor b/2 \rfloor}\frac{1}{(n+2-a-j)(n+2-a-(b+1-j))}\right]\\  & & \mbox{}\times \left[\frac{1}{\lfloor b/2 \rfloor}\sum_{j=1}^{\lfloor b/2 \rfloor}(b - 1 + \kappa_{j-1} + \kappa_{b - j} - \kappa_b)\right].
\end{eqnarray*}

If $b$ is even, then by symmetry we have
$$
\sum_{j=1}^{b/2}(b - 1 + \kappa_{j-1} + \kappa_{b - j} - \kappa_b) =  \frac{1}{2}\sum_{j=1}^b(b - 1 + \kappa_{j-1} + \kappa_{b - j} - \kappa_b) = 0,
$$
and $\Delta_2(n,a,b) \ge 0$, as desired.

If $b$ is odd, then again by symmetry we have
\begin{eqnarray*}
\lefteqn{\sum_{j=1}^{(b-1)/2}(b - 1 + \kappa_{j-1} + \kappa_{b - j} - \kappa_b)}\\
& = & \frac{1}{2}\left[\sum_{j=1}^b(b - 1 + \kappa_{j-1} + \kappa_{b - j} - \kappa_b) - (b -1 + 2 \kappa_{(b-1)/2}-\kappa_b)\right]\\
 & = & -\frac{1}{2}(b -1 + 2 \kappa_{(b-1)/2}-\kappa_b).
\end{eqnarray*}
Hence
\begin{eqnarray*}
\lefteqn{\Delta_2(n,a,b)}\\
& \ge & -[2(n+2-a)-(b+1)]\frac{1}{(b-1)/2}\sum_{j=1}^{(b-1)/2}\frac{1}{(n+2-a-j)(n+2-a-(b+1-j))}\\
& & \mbox{}\times  (b -1 + 2 \kappa_{(b-1)/2}-\kappa_b) + \frac{2}{n+2-a-\frac{b+1}{2}}(b - 1 + 2\kappa_{\frac{b-1}{2}} - \kappa_b)\\
& = & [\kappa_b-2 \kappa_{\frac{b-1}{2}}-(b-1)]\\
& & \mbox{}\times  \Biggl\{\frac{2(n+2-a)-(b+1)}{(b-1)/2}\sum_{j=1}^{(b-1)/2}\frac{1}{(n+2-a-j)(n+2-a-(b+1-j))}\\
& & \qquad \mbox{}-  \frac{4}{2(n+2-a)-(b+1)}\Biggl\}
\end{eqnarray*}
\begin{eqnarray*}
\lefteqn{}\\
& = & [\kappa_b-2 \kappa_{(b-1)/2}-(b-1)][2(n+2-a)-(b+1)]\\
& & \times \Biggl\{ \frac{1}{(b-1)/2}\sum_{j=1}^{(b-1)/2}\frac{1}{(n+2-a-j)(n+2-a-(b+1-j))}\\
& &\quad\mbox{}- \frac{4}{[2(n+2-a)-(b+1)]^2} \Biggl\}\\
& \ge & [\kappa_b-2 \kappa_{(b-1)/2}-(b-1)][2(n+2-a)-(b+1)]\\
& &\times \left\{ \frac{1}{(n+2-a-(b-1)/2)(n+2-a-(b+3)/2)} - \frac{4}{[2(n+2-a)-(b+1)]^2} \right\}
\end{eqnarray*}
\begin{eqnarray*}
\lefteqn{}\\
& = & 4[\kappa_b-2 \kappa_{(b-1)/2}-(b-1)][2(n+2-a)-(b+1)]\\
& & \mbox{}\times \left\{ \frac{1}{[2(n+2-a)-(b-1)][2(n+2-a)-(b+3)]} - \frac{1}{[2(n+2-a)-(b+1)]^2} \right\}\\
& =  & \frac{16[\kappa_b-2 \kappa_{(b-1)/2}-(b-1)][2(n+2-a)-(b+1)]}{[2(n+2-a)-(b-1)][2(n+2-a)-(b+3)][2(n+2-a)-(b+1)]^2}  \\
& \ge & 0,
\end{eqnarray*}
the last two inequalities following from the facts that $2(n+2-a)-(b+3)$ and $[\kappa_b-2 \kappa_{(b-1)/2}-(b-1)]$ are nonnegative---the former due to the fact that $n \ge a+b$, and the latter to~\eqref{orig1} and the identity
\begin{equation}
\label{basicng}
\sum_{j=1}^b(b - 1 + \kappa_{j-1} + \kappa_{b - j} - \kappa_b) = 0.
\end{equation}
\end{proof}

\begl
\label{delt1lmng}
If
\begin{eqnarray*}
 \Delta_1(n,a,b) & := & 2\sum_{j=1}^{b}H_{n+1-j-a}(b - 1 + \kappa_{j-1} + \kappa_{b - j} - \kappa_b),
\end{eqnarray*}
for any nonnegative integers $a$, $b$, and $n$ with $a+b \le n$, then we have
\begin{eqnarray*}
\Delta_1(n,a,b) \le 0.
\end{eqnarray*}
\enl

\begin{proof}
We have
\begin{eqnarray*}
 \Delta_1(n+1,a,b) - \Delta_1(n,a,b) & = & 2\sum_{j=1}^{b}(H_{n+2-j-a} - H_{n+1-j-a})(b - 1 + \kappa_{j-1} + \kappa_{b - j} - \kappa_b)\\
 & = & 2\sum_{j=1}^{b}\frac{1}{n+2-j-a}(b - 1 + \kappa_{j-1} + \kappa_{b - j} - \kappa_b)\\
 & = & \Delta_2(n,a,b).
\end{eqnarray*}
It follows from \refL{delt2lmng} that $\Delta_1(n,a,b)$ is nondecreasing in $n \ge a+b$. From the 
identity~\eqref{basicng} we have
$$
\Delta_1(n,a,b) = 2\sum_{j=1}^b(H_{n+1-j-a}-H_{n})(b - 1 + \kappa_{j-1} + \kappa_{b - j} - \kappa_b).
$$
As a result,
$$
\lim_{n \to \infty}\Delta_1(n,a,b) = 0
$$
follows from the fact that
$$
\lim_{n \to \infty}(H_{n+1-j-a}-H_{n}) = 0
$$
for any $1 \le j \le b$. Thus $\Delta_1(n,a,b) \le 0$ for any $n \ge a+b$, which finishes the proof of the lemma.
\end{proof}

\begl
\label{kaplimgn}
For any nonnegative integer $b$, we have
$$
\lim_{m \to \infty}\sum_{j=1}^{b}\kappa_{m-j}(b - 1 + \kappa_{j-1} + \kappa_{b - j} - \kappa_b) = 0.
$$
\enl

\begin{proof}
Since
\begin{eqnarray*}
\kappa_m 
& = & 2(m+1)H_m - 4m = 2(m+1)\left[\ln m +\gamma+\frac{1}{2m}+O\left(\frac{1}{m^2}\right)\right]-4m\\
& = & 2m\ln m -(4-2\gamma)m+2\ln m +(2\gamma+1)+O\left(\mbox{$\frac{1}{m}$}\right),
\end{eqnarray*}
we have for each fixed~$j$ that
\begin{eqnarray*}
\kappa_{m-j} & = & 2(m-j)\ln(m-j)-(4-2\gamma)(m-j)\\
& & \mbox{}+2\ln(m-j)+(2\gamma+1)+O\left(\mbox{$\frac{1}{m-j}$}\right)\\
& = & 2(m-j)\left[\ln m +\ln\left(1-\sfrac{j}{m}\right)\right]-(4-2\gamma)(m-j)\\
&  & \mbox{}+ 2\left[\ln m +\ln\left(1-\sfrac{j}{m}\right)\right]+(2\gamma+1)+O\left(\sfrac{1}{m-j}\right)\\
& = & 2(m-j)\ln m -2j-(4-2\gamma)(m-j)+2\ln m +(2\gamma+1)+O\left(\sfrac{1}{m}\right)\\
& = & 2m\ln m -(4-2\gamma)m-2(j-1)\ln m +2(1-\gamma)j+(2\gamma+1)+O\left(\sfrac{1}{m}\right).
\end{eqnarray*}
So
\begin{eqnarray*}
\lefteqn{\sum_{j=1}^b\kappa_{m-j}(b - 1 + \kappa_{j-1} + \kappa_{b - j} - \kappa_b)}\\
 & = & \sum_{j=1}^b\left[2m\ln m -(4-2\gamma)m-2(j-1)\ln m +2(1-\gamma)j+(2\gamma+1)+O\left(\sfrac{1}{m}\right)\right]\\
 & & \mbox{}\times  (b - 1 + \kappa_{j-1} + \kappa_{b - j} - \kappa_b)\\
& = & -2[\ln m-(1-\gamma)]\sum_{j=1}^bj(b - 1 + \kappa_{j-1} + \kappa_{b - j} - \kappa_b) + O\left(\sfrac{1}{m}\right),
\end{eqnarray*}
thanks once more to the identity~\eqref{basicng}.

Now
\begin{eqnarray*}
\lefteqn{\sum_{j=1}^{b}j(b - 1 + \kappa_{j-1} + \kappa_{b - j} - \kappa_b)}\\
 & = & \frac{1}{2} \sum_{j=1}^{b}[j+(b+1-j)](b - 1 + \kappa_{j-1} + \kappa_{b - j} - \kappa_b)\\
 & = & \frac{1}{2}(b+1) \sum_{j=1}^{b}(b - 1 + \kappa_{j-1} + \kappa_{b - j} - \kappa_b) = 0
\end{eqnarray*}
by symmetry and~\eqref{basicng}, which finishes the proof of the lemma.
\end{proof}

\begl\label{lam1ng} If
$$
\Lambda_1(a,b) := 2 \sum_{j=1}^{b}H_{j+a} (b - 1 + \kappa_{j-1} + \kappa_{b - j} - \kappa_b)
$$
for any nonnegative integers $a$ and $b$, then
$$
\Lambda_1(a,b) \leq 0.
$$
\enl

\begin{proof}
We have
 \begin{eqnarray*}
\Lambda_1(a,b) & = & 2\sum_{j=1}^{b}H_{b-j+a+1}(b - 1 + \kappa_{j-1} + \kappa_{b - j} - \kappa_b)\\
& = & \Delta_1(b+2a,a,b) \le 0,
\end{eqnarray*}
where the inequality follows from Lemma \ref{delt1lmng}.
\end{proof}

\begl\label{lamng} If
\begin{eqnarray*}
\Lambda(a,b) & := & \sum_{j=1}^{b}\kappa_{j+a-1}(b - 1 + \kappa_{j-1} + \kappa_{b - j} - \kappa_b)\\
& = & \sum_{j=1}^{b}\kappa_{a+b-j}(b - 1 + \kappa_{j-1} + \kappa_{b - j} - \kappa_b)
\end{eqnarray*}
for any nonnegative integers $a$ and $b$, then
$$
\Lambda(a,b) \geq 0.
$$
\enl

\begin{proof}
We have
\begin{eqnarray*}
\Lambda(a+1,b)-\Lambda(a,b) & = & 2\sum_{j=1}^{b}H_{j+a}(b - 1 + \kappa_{j-1} + \kappa_{b - j} - \kappa_b)\\
& = & \Lambda_1(a,b);
\end{eqnarray*}
so, in light of Lemma \ref{lam1ng}, $\Lambda(a,b)$ is nonincreasing in $a$.

We also have from Lemma \ref{kaplimgn} that
$$
\lim_{a \to \infty}\Lambda(a,b) = 0,
$$
which finishes the proof.
\end{proof}

\begl\label{sigmlmng} If
\begin{eqnarray*}
 \Sigma(n,a,b) & := & \sum_{j=a+1}^{a+b}(\kappa_{j-1} + \kappa_{n-j})(b - 1 + \kappa_{j-1-a} + \kappa_{b + a - j} - \kappa_b)\nonumber\\
 & = & \sum_{j=1}^{b}(\kappa_{j+a-1} + \kappa_{n-j-a})(b - 1 + \kappa_{j-1} + \kappa_{b - j} - \kappa_b)
 \end{eqnarray*}
for any nonnegative integers~$a$, $b$, and~$n$ with $a + b \leq n$, then
$$
\Sigma(n,a,b) \ge 0.
$$
\enl

\begin{proof}
We have that
\begin{eqnarray*}
\lefteqn{\Sigma(n+1,a,b) - \Sigma(n,a,b)}\nonumber\\
 & = & \sum_{j=1}^b(\kappa_{n+1-j-a} - \kappa_{n-j-a})(b - 1 + \kappa_{j-1} + \kappa_{b - j} - \kappa_b)\nonumber\\
 & = & 2\sum_{j=1}^bH_{n+1-j-a}(b - 1 + \kappa_{j-1} + \kappa_{b - j} - \kappa_b) = \Delta_1(n,a,b),
 \end{eqnarray*}
which implies from \refL{delt1lmng} that $\Sigma(n, a, b)$ is nonincreasing in $n \geq a + b$. Recall that
\begin{eqnarray*}
\lefteqn{\Sigma(n,a,b)} \\
  & = & \sum_{j=1}^b\kappa_{j+a-1}(b - 1 + \kappa_{j-1} + \kappa_{b - j} - \kappa_b) 
  + \sum_{j=1}^b\kappa_{n-j-a}(b - 1 + \kappa_{j-1} + \kappa_{b - j} - \kappa_b) \\
  & = & \Lambda(a,b) + \sum_{j=1}^b\kappa_{n-j-a}(b - 1 + \kappa_{j-1} + \kappa_{b - j} - \kappa_b);
 \end{eqnarray*}
so, from \refL{kaplimgn} we have
$$
\lim_{ n \to \infty}\Sigma(n,a,b) = \Lambda(a,b).
$$
The result follows from \refL{lamng}.
\end{proof}

\section{The random variables $S_w(t)$, $w \in \Sigma^*$, are nonnegatively correlated}
\label{mainr}

In this section, we first prove the following (in Subsection~\ref{mainr1}) and then complete the proof 
of \refP{mainngthm2} in Subsection~\ref{g1}.  

\begp
\label{Sempng} 
Let $w \in \Sigma^*$. Then the random variables $S_{\emptyset} (t)$ and $S_w(t)$ are nonnegatively correlated.
\enp

\subsection{The random variables $S_{\emptyset} (t)$ and $S_w(t)$ for any $w \in \Sigma^*$ are nonnegatively correlated}
\label{mainr1}

In this Subsection~\ref{mainr1} we prove \refP{Sempng}, which states that
\begin{equation}
\label{cov1}
\Cov(S_{\emptyset} (t), S_w(t)) \geq 0\,\,\text{for any}\,\, w \in \Sigma^*,
\end{equation}
with the understanding that $S_{\emptyset} (t) = K(t) = K_{N(t)}$.

\begin{proof}[\textbf{Proof of \refP{Sempng}}]
We have
\begin{equation}
\label{cov2}
\Cov(S_{\emptyset} (t), S_w(t)) = \Cov(K(t), S_w(t)) = T_w(t) + V_w(t)
\end{equation}
where
\begin{equation}
\label{Twtng}
T_w(t) := \Cov(\EE[K(t)\,|\,N(t)], \EE[S_w(t)\,|\,N(t)])
\end{equation}
and
\begin{equation}
\label{Vwtng}
V_w(t) := \EE\,\Cov(K(t),S_w(t)\,|\,N(t)).
\end{equation}
But Propositions \ref{Twposthm}--\ref{Vwposthm} will demonstrate that the expressions $T_w(t)$ and $V_w(t)$ are each nonnegative. 
\end{proof}

\subsubsection{Nonnegativity of $T_w(t)$}\label{Twposng}

Here we prove the following result.

\begp
\label{Twposthm}
The expression $T_w(t)$ defined in~\eqref{Twtng} is nonnegative.
\enp

\begin{proof}
We have
$$
\EE[K(t)\,|\,N(t) = n] = \kappa_n := \EE\,K_n,
$$
which is increasing with $n$; and
$$
\EE[S_w(t)\,|\,N(t) = n] = \sum_{j=0}^n \binom{n}{j} p_w^j (1 - p_w)^{n - j} \kappa_j
$$
is also increasing, following from the fact that the Binomial$(n, p_w)$ distributions increase stochastically with~$n$.

By ``Chebyshev's other inequality"~\cite{MR789242}, we can conclude that
$$
\Cov(\EE[K(t)\,|\,N(t)], \EE[S_w(t)\,|\,N(t)]) \geq 0,
$$
which finishes the proof of the proposition.
\end{proof}

\subsubsection{Nonnegativity of $V_w(t)$}\label{Vwposng}

In this subsection we prove the following proposition, thereby completing the proof of \refP{Sempng}.

\begp
\label{Vwposthm}
The expression $V_w(t)$ defined in \eqref{Vwtng} is nonnegative.
\enp
\noindent
This will be accomplished using the next two propositions, Propositions \ref{psing} and \ref{ching}.

\begp
\label{psing} 
If
\begin{eqnarray}\label{psidefng}
\lefteqn{\psi(n,a,b)} \nonumber \\
& := & n^{-1}\!\!\!\!\!\!\sum_{a < j \le a+ b}(n-1 + \kappa_{j-1} + \kappa_{n-j} - \kappa_n)
(b - 1 + \kappa_{j-1-a} + \kappa_{b + a - j} - \kappa_{b})
\end{eqnarray}
for any nonnegative integers $a$, $b$, and $n$ with $a+b \le n$, then
$$
\psi(n,a,b) \geq 0.
$$
\enp

\begin{proof}
We have
\begin{eqnarray*}
\psi(n,a,b) & = & n^{-1} \sum_{a < j \le a + b}(\kappa_{j-1} + \kappa_{n-j})(b- 1 + \kappa_{j-1-a} + \kappa_{b + a - j} - \kappa_{b})\nonumber\\
 & & \quad \mbox{} + n^{-1} (n - 1 - \kappa_n) \sum_{a < j \le a + b}(b - 1 + \kappa_{j-1-a} + \kappa_{b + a - j} - \kappa_{b})\nonumber\\
 & = & n^{-1} \sum_{a < j \le a + b}(\kappa_{j-1} + \kappa_{n-j})(b - 1 + \kappa_{j-1-a} + \kappa_{b + a - j} - \kappa_{b})\nonumber\\
 & = & n^{-1} \Sigma(n,a,b) \geq 0,
\end{eqnarray*}
where the second equality follows from the fact that
\begin{equation}
\label{orig}
\kappa_{d} = \sfrac{1}{d}\sum_{j=1}^{d}(d - 1 + \kappa_{j-1} + \kappa_{d - j}) = \sfrac{1}{d} \sum_{j=e+1}^{e + d}(d - 1 + \kappa_{j-1-e} + \kappa_{e + d - j})
\end{equation}
for any two integers $e$ and $d \ge 1$, and the inequality from Lemma \ref{sigmlmng}.
\end{proof}

\begin{definition}
\label{Snwng}
\emph{
Let $w \in \sum^*$, and let $n$ be any nonnegative integer. We define $S_{n,w}$ to be the number of key comparisons between those keys (from among the $n$ first to arrive) with prefix $w$.
}
\end{definition}

\begin{definition}
\label{Nnwng}
\emph{
For any $w \in \sum^*$, and nonnegative integer $n$, we define $N_{n,w}$ to be the number of keys (from among the $n$ first to arrive) with prefix $w$, and
$$
N_{n,w^-} := \sum_{\genfrac{}{}{0pt}{}{w'\in \Sigma^{|w|}:}{w'<w}}N_{n,w'}.
$$
}
\end{definition}

\begp
\label{ching} 
For any nonnegative integers $a$, $b$, and $n$ with $a+b \le n$, we have
\begin{equation}
\label{chieqn}
\Cov(K_n,S_{n,w}\,|\,N_{n,w} = b,N_{n,w^-} = a) \geq 0.
\end{equation}
\enp

\begin{proof}
We will prove the proposition by strong induction on~$n$.
For that, we further condition on $J_n :=$\ (the rank of the root key among the first $n$ keys). Applying the law of total covariance 
(namely, covariance equals the sum of expectation of conditional covariance and covariance of conditional expectations) 
to the conditional covariance in question, we find
\begin{eqnarray}
\lefteqn{\hspace{.4in}\Cov(K_n, S_{n,w}\,|\,N_{n,w} = b,N_{n,w^-} = a)} \label{cov4} \\
 & = & \sum_{j=1}^n \PP[J_n = j|  N_{n,w} = b,N_{n,w^-} = a] \nonumber\\
& & \qquad \times (\EE[K_n |  N_{n,w} = b,N_{n,w^-} = a,J_n = j]-\EE[K_n |  N_{n,w} = b,N_{n,w^-} = a])\nonumber\\
& & \qquad \times  (\EE[S_{n,w} |  N_{n,w} = b,N_{n,w^-} = a,J_n = j]-\EE[S_{n,w} |  N_{n,w} = b,N_{n,w^-} = a]) \nonumber\\
&  & \mbox{}+ \sum_{j=1}^n \PP[J_n = j|  N_{n,w} = b,N_{n,w^-} = a] \nonumber \\
&  & \qquad \times \Cov(K_n,S_{n,w}|  N_{n,w} = b,N_{n,w^-} = a,J_n = j).\nonumber
\end{eqnarray}

In preparation for handling~\eqref{cov4}, we begin with three observations, mainly concerning the first of the two terms on the right in~\eqref{cov4}.
\medskip

{\bf (i)}~$(K_n, J_n)$ and $( N_{n,w}, N_{n,w^-})$ are independent, so for any $j = 1, \ldots, n$ and any nonnegative integers~$a$ and~$b$, we have
$$
\PP[J_n = j\,|\,N_{n,w} = b, N_{n,w^-} = a]  =  \PP[J_n = j] = \sfrac{1}{n},
$$
and
$$
\EE[K_n\,|\,N_{n,w} = b, N_{n,w^-} = a, J_n = j] = \EE[K_n\,|\,J_n = j] = n - 1 + \kappa_{j - 1} + \kappa_{n - j},
$$
and
$$
\EE[K_n\,|\,N_{n,w} = b, N_{n,w^-} = a] = \EE\,K_n = \kappa_n.
$$
Also
$$
\EE[S_{n,w}\,|\,N_{n,w} = b, N_{n,w^-} = a] = \kappa_b.
$$

Keep in mind in the observations to follow that~$a$ is the value of $N_{n,w^-}$, that~$b$ is the value of $N_{n,w}$, and that~$j$ is the value of $J_n$.
\medskip

{\bf (ii)}~If
$a < j \leq a + b$, which happens in the case that the root key has its prefix of length $|w|$ equal to~$w$, then there are $j - 1 - a$ keys among the $j - 1$ that fall to the left of the pivot key that have~$w$ as their prefix of
length $|w|$, and  $b + a - j$ keys among the $n - j$ that fall to the right of the pivot key that have~$w$ as their prefix of
length $|w|$. So
$$
\mathcal{L}(S_{n,w}\,|\,N_{n,w} = b,N_{n,w^-} = a,J_n = j) = \mathcal{L}(b - 1 + D'_{j-1,j-1-a} + D''_{n-j,b + a - j})
$$
where $D'_{j-1,j-1-a}$ and $D''_{n-j,b + a - j}$ are independent, and
$$
\mathcal{L}(D'_{j-1,j-1-a}) = \mathcal{L}(S_{j-1,w}|  N_{j-1,w} = j-1-a,N_{j-1,w^-} = a)
= \mathcal{L}(K_{j-1-a}),
$$
and similarly 
$$
\mathcal{L}(D''_{n-j,b + a - j}) = \mathcal{L}(K_{b + a - j});
$$
hence
$$
\EE[S_{n,w}\,|\,N_{n,w} = b, N_{n,w^-} = a, J_n = j] = b - 1 + \kappa_{j-1-a} + \kappa_{b+a - j}.
$$

{\bf (iii)}~If $j \leq a$ or $a + b < j$, which happens if the root key has its prefix of length $|w|$ different from~$w$, then all of the keys that have~$w$ as their prefix of length $|w|$ fall on the same side of the pivot key.  So
$$
\mathcal{L}(S_{n,w}\,|\,N_{n,w} = b, N_{n,w^-} = a, J_n = j) = \mathcal{L}(K_b)
$$
and
$$
\EE[S_{n,w}\,|\,N_{n,w} = b, N_{n,w^-} = a, J_n = j] = \kappa_b.
$$

Equation~\eqref{cov4} now yields
\begin{eqnarray*}
\lefteqn{\Cov(K_n, S_{n,w}\,|\,N_{n,w} = b, N_{n,w^-} = a)} \\
 & = & \frac{1}{n}\biggl\{\sum_{a < j \le a+ b}(n-1 + \kappa_{j-1} + \kappa_{n-j} - \kappa_n)
(b - 1 + \kappa_{j-1-a} + \kappa_{b + a - j} - \kappa_b)\\
 &  & \qquad \mbox{}+\sum_{1 \le j \le a}(n-1 + \kappa_{j-1} + \kappa_{n-j} - \kappa_n)(\kappa_b - \kappa_b) \\
 &  & \qquad \mbox{}+\sum_{a+b <j \le n}(n-1 + \kappa_{j-1} + \kappa_{n-j} - \kappa_n)(\kappa_b - \kappa_b)\biggl\}\\
&  & \quad \mbox{}+ \frac{1}{n}\sum_{j=1}^n\Cov(K_n, S_{n,w}\,|\,N_{n,w} = b,N_{n,w^-} = a,J_n = j)\\
& = & \frac{1}{n}\sum_{a < j \le a+ b}(n-1 + \kappa_{j-1} + \kappa_{n-j} - \kappa_n)(b - 1 + \kappa_{j-1-a} + \kappa_{b + a - j} - \kappa_b)\\
&  & \quad \mbox{}+ \frac{1}{n}\sum_{j=1}^n\Cov(K_n, S_{n,w}\,|\,N_{n,w} = b,N_{n,w^-} = a,J_n = j)\\
& = & \psi(n,a,b) + \frac{1}{n}\sum_{j=1}^n\Cov(K_n, S_{n,w}\,|\,N_{n,w} = b,N_{n,w^-} = a,J_n = j)\\
& \geq & \frac{1}{n}\sum_{j=1}^n\Cov(K_n, S_{n,w}\,|\,N_{n,w} = b,N_{n,w^-} = a,J_n = j),
\end{eqnarray*}
where the last equality follows from~\eqref{psidefng}, and the inequality from \refP{psing}. So, to prove 
that~\eqref{chieqn} holds, we only need to prove that
\begin{equation}
\label{KSjeq}
\Cov(K_n, S_{n,w}\,|\,N_{n,w} = b, N_{n,w^-}=a, J_n = j) \geq 0\,\,\text{for any}\,\,1 \le j \le n.
\end{equation}

First note that if $n=1$, then $K_n \equiv 0$ and hence~\eqref{KSjeq} holds.

Now let's assume that~\eqref{chieqn} holds for any natural number smaller than a given natural number~$n$. Then:
\medskip

\noindent
{\sc Case A.}\ \ 
If $a < j \leq a + b$ \,then there are $j - 1 - a$ keys among the $j - 1$ that fall to the left of the pivot key that have  their prefix of length $|w|$ equal to~$w$, and  $b + a - j$ keys among the $n - j$ that fall to the right of the pivot key that have  their prefix of length $|w|$ equal to~$w$.  So
\begin{eqnarray*}
\lefteqn{\mathcal{L}(K_n,S_{n,w}| N_{n,w} = b,N_{n,w^-}=a,J_n = j)}\\
 & = & \mathcal{L}(n-1+K'_{j-1}+K''_{n-j},b - 1 + D'_{j-1,j-1-a} + D''_{n-j,b + a - j})
 \end{eqnarray*}
 where
$$
\mathcal{L}(K'_{j-1},D'_{j-1,j-1-a}) = \mathcal{L}(K_{j-1},S_{j-1,w}| N_{j-1,w} = j-1-a,N_{j-1,w^-} = a)
$$
and 
$$
\mathcal{L}(K''_{n-j},D''_{n-j,b + a - j}) = \mathcal{L}(K_{n-j},S_{n-j,w}| N_{n-j,w} =b+a-j,N_{n-j,w^-}=0)
$$
and also 
$$
\mbox{$(K'_{j-1},D'_{j-1,j-1-a})$ and $(K''_{n-j},D''_{n-j,b + a - j})$ are independent}.
$$

In this case, therefore,
\begin{eqnarray*}
\lefteqn{\Cov(K_n,S_{n,w}\,|\,N_{n,w} = b,N_{n,w^-}=a,J_n = j)}\\
 & = & \Cov(n-1+K'_{j-1}+K''_{n-j}, b - 1 + D'_{j-1,j-1-a} + D''_{n-j,b+ a - j})\\
& = & \Cov(K'_{j-1}, D'_{j-1,j-1-a}) + \Cov(K''_{n-j}, D''_{n-j,b + a - j})\\
& = & \Cov(K_{j-1}, S_{j-1,w}\,|\,N_{j-1,w} = j-1-a, N_{j-1,w^-}=a)\\
&  & \quad \mbox{}+ \Cov(K_{n-j}, S_{n-j,w}\,|\,N_{n-j,w} =b+a-j, N_{n-j,w^-}=0) \geq 0
\end{eqnarray*}
by strong induction, since $j-1 < n$ and $n - j < n$.
\medskip

\noindent
{\sc Case~B.}\ \ 
If $j \leq a$, which happens if the keys that have $w$ as their prefix of length $|w|$ all fall to the right of the pivot key, then
$$
\mathcal{L}(K_n, S_{n,w}\,|\,N_{n,w} = b, N_{n,w^-}=a, J_n = j) = \mathcal{L}(n-1+K'_{j-1}+K''_{n-j}, D''_{n-j,b})
$$
where
$$
\mathcal{L}(K'_{j-1})=\mathcal{L}(K_{j-1})
$$
and  
$$
\mathcal{L}(K''_{n-j},D''_{n-j,b}) = \mathcal{L}(K_{n-j},S_{n-j,w}\,|\,N_{n-j,w} = b,N_{n-j,w^-} = a-j)
$$
and also
$$
\mbox{$K'_{j - 1}$ and $(K''_{n-j}, D''_{n-j,b})$ are independent}.
$$

In this case, therefore,
\begin{eqnarray*}
\lefteqn{\Cov(K_n, S_{n,w}\,|\,N_{n,w} = b, N_{n,w^-}=a, J_n = j)}\\
& = & \Cov(n-1+K'_{j-1}+K''_{n-j}, D''_{n-j,b}) = \Cov(K''_{n-j}, D''_{n-j,b})\\
& = & \Cov(K_{n-j}, S_{n-j,w}\,|\,N_{n-j,w} = b, N_{n-j,w^-} = a-j) \geq 0
\end{eqnarray*}
by strong induction, since $n - j < n$.
\medskip

\noindent
{\sc Case~C.}\ \ 
If $a + b < j$, which happens if the keys that have~$w$ as their prefix of length $|w|$ all fall to the left of the pivot key, then
$$
\mathcal{L}(K_n,S_{n,w}| N_{n,w} = b,N_{n,w^-}=a,J_n = j) = \mathcal{L}(n-1+K'_{j-1}+K''_{n-j}, D'_{j-1,b})
$$
where
$$
\mathcal{L}(K'_{j-1},D'_{j-1,b})=\mathcal{L}(K_{j-1},S_{j-1,w}| N_{j-1,w} = b,N_{j-1,w^-} = a)
$$
and 
$$
\mathcal{L}(K''_{n-j})=\mathcal{L}(K_{n-j})
$$
and also 
$$
\mbox{$(K'_{j-1}, D'_{j-1,b})$ and $K''_{n-j}$ are independent}.
$$

In this case, therefore,
\begin{eqnarray*}
\lefteqn{\Cov(K_n, S_{n,w}\,|\,N_{n,w} = b, N_{n,w^-}=a, J_n = j)}\\
& = & \Cov(n-1+K'_{j-1}+K''_{n-j}, D'_{j-1,b}) = \Cov(K'_{j-1}, D'_{j-1,b})\\
& = & \Cov(K_{j-1},S_{j-1,w}\,|\,N_{j-1,w} = b, N_{j-1,w^-} = a) \geq 0
\end{eqnarray*}
by strong induction, since $j - 1 < n$.
\medskip

In all three cases~\eqref{KSjeq} holds, which concludes the proof of the proposition.
\end{proof}

\begin{proof}[\textbf{Proof of \refP{Vwposthm}}]
To prove \refP{Vwposthm}, which asserts that
$$
\EE\,\Cov(K(t), S_w(t)\,|\,N(t)) \geq 0,
$$
it's enough to show that
$$
\Cov(K(t), S_w(t)\,|\,N(t)=n) \geq 0\,\,\text{ for all}\,\, n=0, 1, 2, \ldots.
$$
But
$$
\Cov(K(t),S_w(t)\,|\,N(t)=n) = \Cov(K_n, S_{n,w}),
$$
and conditioning on $N_{n,w}$ and $N_{n,w^-}$ we have
\begin{eqnarray*}
\Cov(K_n, S_{n,w}) & = & \Cov(\EE[K_n\,|\,N_{n,w}, N_{n,w^-}], \EE[S_{n,w}\,|\,N_{n,w}, N_{n,w^-}]) \\
 & & \qquad \mbox{}+ \EE\,\Cov(K_n, S_{n,w}\,|\,N_{n,w}, N_{n,w^-}).
\end{eqnarray*}
Knowing that $K_n $ and $(N_{n,w}, N_{n,w^-})$ are independent, we have
$$
\Cov(\EE[K_n\,|\,N_{n,w}, N_{n,w^-}], \EE[S_{n,w}\,|\,N_{n,w}, N_{n,w^-}]) = \Cov(\kappa_n, \kappa_{N_{n,w}}) = 0.
$$
We have now reduced to proving
$$
\EE\,\Cov(K_n, S_{n,w}\,|\,N_{n,w}, N_{n,w^-}) \geq 0,
$$
which is achieved by~\refP{ching}.
\end{proof}

\subsection{The general case}\label{g1}

\begin{proof}[\textbf{Proof of Proposition \ref{mainngthm2}}]

Let $w$ and $w'$ be in $\Sigma^*$. 
On the one hand, if the prefixes~$w$ and $w'$ are inconsistent in the sense that no word has both~$w$ and $w'$ as prefixes (for example, if $w = 01$ and $w' = 1$), then $S_w(t)$ and $S_{w'}(t)$ are independent and therefore uncorrelated.  On the other hand, if~$w$ and $w'$ are not inconsistent, then either $w'$ is a prefix of~$w$ or~$w$ is a prefix of~$w'$ (or both, which is precisely the case $w = w'$).  Let's assume without loss of generality that $w'$ is a prefix of~$w$; then $w = w' w''$, the concatenation of $w'$ with another prefix $w''$.  Having begun with a probabilistic source~$\mu$, consider the source~$\mu'$ obtained by conditioning on prefix~$w'$, and use notation $S'$ for 
symbol-count variables for source $\mu'$ just as~$S$ is used for source~$\mu$.  [Observe that $\mu'$, like~$\mu$, satisfies the condition~\eqref{cond}.]  Then
$$
\mathcal{L}(S_{w'}(t), S_w(t)) = \mathcal{L}(S'_\emptyset(p_{w'} t), S'_{w''}(p_{w'} t)).
$$
The result follows from \refP{Sempng}.
\end{proof}

\bibliographystyle{plain}
\bibliography{references}
\end{document}